\documentclass[12pt]{amsart}
\usepackage{amsmath, amsthm, amscd, amsfonts}

\newtheorem{theorem}{Theorem}[section]
\newtheorem{lemma}[theorem]{Lemma}

\theoremstyle{definition}
\newtheorem{definition}[theorem]{Definition}
\newtheorem{example}[theorem]{Example}

\theoremstyle{remark}
\newtheorem{remark}[theorem]{Remark}
\numberwithin{equation}{section}
\newfont{\kh}{msbm10}
\newcommand{\R}{\mbox{\kh R}}
\newcommand{\C}{\mbox{\kh C}}
\begin{document}
\title{Dynamical Systems on Hilbert $C^*$-Modules}
\author{Gh. Abbaspour Tabadkan}
\address{Gholamreza Abbaspour Tabadkan, \newline Department of Mathematics, Ferdowsi University, P. O. Box 1159, Mashhad 91775, Iran}
\email{tabadkan@math.um.ac.ir}
\author{M. S. Moslehian}
\address{Mohammad Sal Moslehian, \newline Department of Mathematics, Ferdowsi University, P. O. Box 1159, Mashhad 91775, Iran}
\email{moslehian@ferdowsi.um.ac.ir}
\author{A. Niknam}
\address{Assadollah Niknam, \newline Department of Mathematics, Ferdowsi University, P. O. Box 1159, Mashhad 91775, Iran}
\email{niknam@math.um.ac.ir} 
\begin{abstract}
We investigate the generalized derivations and show that every
generalized derivation on a simple Hilbert $C^*$-module either is
closable or has a dense range. We also describe dynamical systems
on a full Hilbert $C^*$-module ${\mathcal M}$ over a $C^*$-algebra
${\mathcal A}$ as a one-parameter group of unitaries on ${\mathcal
M}$ and prove that if $\alpha: \R\to U({\mathcal M})$ is a
dynamical system, where $U({\mathcal M})$ denotes the set of all
unitary operator on ${\mathcal M}$, then we can correspond a
$C^*$-dynamical system $\alpha^{'}$ on ${\mathcal A}$ such that if
$\delta$ and $d$ are the infinitesimal generators of $\alpha$ and
$\alpha^{'}$ respectively, then $\delta$ is a $d$-derivation.
\end{abstract}
\maketitle

\section{Introduction.}

A Hilbert $C^*$-module over a $C^*$-algebra ${\mathcal A}$ is an
algebraic left ${\mathcal A}$-module ${\mathcal M}$ equipped with
an ${\mathcal A}$-valued inner product $\langle \cdot,
\cdot\rangle $ which is ${\mathcal A}$-linear in the first and
conjugate linear in the second variable such that ${\mathcal M}$
is a Banach space with respect to the norm $\| x \|=\|\langle
x,x\rangle \|^{1/2}$. The Hilbert module ${\mathcal M}$ is called
full if the closed linear span $\langle {\mathcal M},{\mathcal
M}\rangle $ of all elements of the form $\langle x,y\rangle
~~(x,y\in {\mathcal M})$ is equal to ${\mathcal A}$.

Hilbert $C^*$-modules are first introduced and investigated by I.
Kaplansky \cite{KAP}, M. Rieffel \cite{RIE} and W. Paschke
\cite{PAS}. They are a generalization of Hilbert spaces, but there
are some differences between these two classes. For example, each
operator on a Hilbert space has an adjoint, but a bounded
${\mathcal A}$-module map on a Hilbert ${\mathcal A}$-module is
not adjointable in general \cite{LAN}.

In this paper, we investigate the generalized derivations. This
notion was first appeared in the context of operator algebras
\cite{MAT}. Later, it was introduced in the framework of pure
algebra \cite{HVA}. We shall show that every generalized
derivation $\delta$ on a simple Hilbert $C^*$-module ${\mathcal
M}$ either is closable or has a dense range in ${\mathcal M}$.
This is a generalization of a result of A. Niknam \cite{NIK1}.

We also describe dynamical systems on a full Hilbert $C^*$-module
${\mathcal M}$ over a $C^*$-algebra ${\mathcal A}$ as a
one-parameter group of unitaries on ${\mathcal M}$ and prove that
if $\alpha: \R\to U({\mathcal M})$ is a dynamical system, then we
can correspond a $C^*$-dynamical system $\alpha^{'}$ on ${\mathcal
A}$ such that if $\delta$ and $d$ are the infinitesimal generators
of $\alpha$ and $\alpha^{'}$ respectively, then $\delta$ is a
$d$-generalized derivation.

The reader is referred to \cite{LAN} for more details on Hilbert
$C^*$-modules and to \cite{SAK} for more information on
$C^*$-dynamical systems.

\section{Preliminaries.}

Throughout this section ${\mathcal M}$ and ${\mathcal N}$ are
assumed to be Hilbert modules over $C^*$-algebras ${\mathcal A}$
and ${\mathcal B}$ respectively, and $\varphi:{\mathcal A} \to
{\mathcal B}$ is a morphism of $C^*$-algebras.

A map $\Phi: {\mathcal M}\to {\mathcal N}$ is said to be a
$\varphi$-morphism of Hilbert $C^*$-modules if
$$\langle \Phi(x),\Phi(y)\rangle =\varphi(\langle x,y\rangle ) ~~(x,y \in {\mathcal M})$$
Applying the polarization we can immediately conclude that $\Phi$
is a $\varphi$-morphism if and only if $\langle
\Phi(x),\Phi(x)\rangle =\varphi(\langle x,x\rangle )$ for all $x
\in {\mathcal M}$. Each $\varphi$-morphism is necessarily a linear
operator since for any $x, y, z\in {\mathcal M}$ and $\alpha \in
\C$ we have
\begin{eqnarray*}
\langle \Phi(\alpha x+y),\Phi(z)\rangle &=&\alpha \varphi(\langle x,z\rangle )+\varphi(\langle y,z\rangle )\\
&=&\alpha\langle \Phi(x),\Phi(z)\rangle +\langle
\Phi(y),\Phi(z)\rangle
\end{eqnarray*}
and hence $\langle \Phi(\alpha
x+y)-\alpha\Phi(x)-\Phi(y),\Phi(z)\rangle  = 0$ for all $z\in
{\mathcal M}$. If we replace $z$ with $x, y$ and $\alpha x+y$ then
we can infer that
\[\langle \Phi(\alpha x+y)-\alpha\Phi(x)-\Phi(y),\Phi(\alpha x+y)-\alpha\Phi(x)-\Phi(y)\rangle =0\]
Thus $\Phi(\alpha x+y)-\alpha\Phi(x)-\Phi(y)=0$ and therefore
$\Phi(\alpha x+y)=\alpha\Phi(x)+\Phi(y)$. Similarly every
$\varphi$-morphism is necessarily a module map in the sense that
$\Phi(ax)=\varphi(a)\Phi(x)$ for all $x\in {\mathcal M}$  and for
all $a \in {\mathcal A}$.

Further, let $\psi:{\mathcal B} \to C$ be a morphism of
$C^*$-algebras and let ${\mathcal Q}$ be a Hilbert $C^*$-module
over $C$. If $\Phi: {\mathcal M}\to {\mathcal N}$ is a
$\varphi$-morphism and $\Psi: {\mathcal N}\to {\mathcal Q}$ is a
$\psi$-morphism, then obviously $\Psi\Phi:{\mathcal M}\to
{\mathcal Q}$ is a $\psi\varphi$-morphism of Hilbert
$C^*$-modules.

Following \cite{B-G} we call a map $\Phi:{\mathcal M}\to {\mathcal
N}$ a unitary operator if there exists an injective morphism of
$C^*$-algebras $\varphi:{\mathcal A}\to {\mathcal B}$ such that
$\Phi$ is a surjection $\varphi$-morphism.

\begin{remark} {\bf (i)} If $\Phi: {\mathcal M}\to {\mathcal N}$ be a $\varphi$-morphism of Hilbert
$C^*$-modules and $\varphi$ is injective then $\Phi$ is an
isometry; cf. \cite{B-G}. Thus each unitary operator of Hilbert
$C^*$-modules is necessarily an isometry.

{\bf (ii)} It is known \cite{B-G} that if $\Phi: {\mathcal M}\to
{\mathcal N}$ is a $\varphi$-morphism of Hilbert $C^*$-modules
then $Im\Phi$ is a closed subspace of ${\mathcal N}$ and a Hilbert
$C^*$-module over the $C^*$-algebra $Im\varphi\subset {\mathcal
B}$ such that $\langle Im\Phi,Im\Phi\rangle =\varphi(\langle
{\mathcal M},{\mathcal M}\rangle )$. Moreover if $\Phi $ is
surjective, and if ${\mathcal N}$ is a full ${\mathcal B}$-module,
then $\varphi$ is also surjective. Thus if $\Phi$ is a unitary
then it is surjective and hence $\langle {\mathcal N},{\mathcal
N}\rangle ~=\varphi(\langle {\mathcal M},{\mathcal M}\rangle
)\simeq~\langle {\mathcal M},{\mathcal M}\rangle $. Moreover, if
${\mathcal N}$ is a full Hilbert ${\mathcal B}$-module, then
$\varphi$ is surjective and so it is an isomorphism of
$C^*$-algebras.
\end{remark}

\begin{example} Let ${\mathcal H}$ be a Hilbert space. Then ${\mathcal H}$ can be
regarded as a $K({\mathcal H})$-module via $T \cdot x=T(x)~~(T\in
K({\mathcal H}), x\in {\mathcal H})$. If we define a $K({\mathcal
H})$-inner product on ${\mathcal H}$ via $\langle x,y\rangle
=x\otimes y~~x,y\in {\mathcal H} $, where $x\otimes y(z)=(z,y)x ~
z\in {\mathcal H}$ and $(\cdot, \cdot)$ denotes the complex inner
product on ${\mathcal H}$, then ${\mathcal H}$ can be regarded as
a Hilbert $C^*$-module over $K({\mathcal H})$. In this case,
$U:{\mathcal H}\to {\mathcal H}$ is a unitary as an operator on
Hilbert space ${\mathcal H}$ if and only if $U$ is a unitary
operator on Hilbert $K({\mathcal H})$-module ${\mathcal H}$ in the
above sense. This is a consequence of the facts that $\langle
U(x),U(y)\rangle =U(x)\otimes U(y)=U(x\otimes y)U^{*}$, and that $
AdU:K({\mathcal H}) \to K({\mathcal H})$ defined by
$AdU(V)=UVU^{*}, V\in K({\mathcal H})$ is a $*$-isomorphism and
each $*$-isomorphism on $K({\mathcal H})$ is of this form; cf.
\cite{MUR}.
\end{example}

We denote by $U({\mathcal M})$ the group of all unitary operators
of ${\mathcal M}$ onto ${\mathcal M}$. If ${\mathcal M}$ is full
and $\alpha:{\mathcal M}\to {\mathcal M}$ is a unitary operator
then by Remark 2.1.ii there is a $*$-isomorphism
$\alpha^{'}:{\mathcal A}\to {\mathcal A}$ such that $\alpha$ is a
$\alpha^{'}$-morphism.

We end this section with the following useful lemma which can be
found in \cite{A-N} and \cite{MOS}.

\begin{lemma} Let ${\mathcal M}$ be a full Hilbert module over the
$C^*$-algebra ${\mathcal A}$ and let $a\in {\mathcal A}$. Then
$a=0$ if and only if $ax=0$ for all $x\in {\mathcal M}$.
\end{lemma}

\section {Generalized Derivation}

This section is devoted to study of generalized derivations. Our
aim is to show that every generalized derivation $\delta$ on a
simple Hilbert $C^*$-module either is closable or has a dense
range.

\begin{definition} Let ${\mathcal M}$ be full Hilbert ${\mathcal A}$-module. A
linear map $\delta:D(\delta)\subseteq {\mathcal M}\to {\mathcal
M}$, where $D(\delta)$ is a dense subspace of ${\mathcal M}$, is
called a generalized derivation if there exists a mapping
$d:D(d)\to {\mathcal A}$, where $D(d)$ is a dense subalgebra of
${\mathcal A}$ such that $D(\delta)$ is an algebraic left
$D(d)$-module, and $\delta(ax)=a\delta(x)+d(a)x$ for all $x\in
D(\delta)$ and all $a\in D(d)$.
\end{definition}

In this case $d$ must be a derivation since for any $a,b\in D(d)$
and $x \in D(\delta)$ we have
\[\delta(abx)=ab\delta(x)+d(ab)x.\]
On the other hand
\[\delta(abx)=\delta(a(bx))=a\delta(bx)+d(a)bx=ab\delta(x)+ad(b)x+d(a)bx,\]
whence
\[(d(ab)-(ad(b)+d(a)b)x=0\]
for all $x\in D(\delta)$. Thus by Lemma 2.3 we obtain
$d(ab)=ad(b)+d(a)b$ since $D(\delta)$ is dense in ${\mathcal M}$.

Similarly we can show that $d$ is linear so $d:D(d)\subseteq
{\mathcal A}\to {\mathcal A}$ is a derivation. We call $\delta$ a
$d$-derivation.

Denote by $GDer({\mathcal M})$ the set of all generalized
derivations on ${\mathcal M}$. Then it is easy to see that
$GDer({\mathcal M})$ is a linear space. In fact if $\delta_{1},
\delta_{2}\in GDer({\mathcal M})$, $\delta_{1}$ is a
$d_{1}$-derivation and $\delta_{2}$ is a $d_{2}$-derivation, then
$\alpha \delta_{1}+\beta \delta_{2}$ is a $\alpha d_{1}+\beta
d_{2}$-derivation and so $\alpha \delta_{1}+\beta \delta_{2} \in
GDer({\mathcal M})$. Also the Lie product $[\delta_{1} ,
\delta_{2}]$ is $[d_{1},d_{2}]$-derivation and so
$[\delta_{1},\delta_{2}]\in GDer({\mathcal M}).$

Now we show a similar result as in \cite{NIK1} for generalized
derivations:

\begin{theorem} Let ${\mathcal M}$ be a simple full Hilbert
$C^*$-module in the sense that it has no trivial left ${\mathcal
A}$-submodule and let $\delta:D(\delta)\subseteq {\mathcal M}\to
{\mathcal M}$ be a $d$-derivation. Then either $\delta$ is
closable or the range of $\delta$ is dense in ${\mathcal M}$.
\end{theorem}

\begin{proof} Let $\sigma(\delta)$ be the separating space of
$\delta$ that is
\[\sigma(\delta)=\{x\in {\mathcal M}  ;{\rm there~ is~ a~
sequence} ~(x_{n})~ {\rm in}~ D(\delta)~ {\rm with} ~x_{n}\to 0,
\delta(x_{n})\to x \}\] It is obvious that $\sigma(\delta)$  is a
closed subspace of ${\mathcal M}$. We show that $\sigma(\delta)$
is a left submodule of ${\mathcal M}$. Let $a\in {\mathcal A},
x\in \sigma(\delta)$ thus there exists a sequence
$(x_{n})\subseteq D(\delta)$ such that $x_{n}\to 0$ and
$\delta(x_{n})\to x $ so we have $ ax_{n}\to 0$ and
$\delta(ax_{n})=a\delta(x_{n})+d(a)x_{n}\to ax$\\Thus $ax\in
\sigma(\delta)$. By the hypothesis $\sigma(\delta)=\{0\}$ or
$\sigma(\delta)={\mathcal M}$. Therefore $\delta$ is closable or
range $\delta$ is dense in ${\mathcal M}$.
\end{proof}

\section {Dynamical Systems On Full Hilbert $C^*$-Modules}

We start this section with a basic definition.

\begin{definition} Let ${\mathcal M}$ be a full Hilbert ${\mathcal A}$-module. A map
$\alpha$ from the real line $\R$ to $U({\mathcal M})$ which maps
$t$ to $\alpha_{t}$ is said to be a one-parameter group of
unitaries if

(i) $\alpha_{0}=I$

(ii) $\alpha_{t+s}=\alpha_{t}\alpha_{s}~~(t,s \in \R)$

Further, $\alpha$ is said to be a strongly continuous
one-parameter group of unitaries if, in addition, $\lim_{t\to\ 0}
\alpha_{t}(x)= x$ in the norm of ${\mathcal M}$ for all $x\in
{\mathcal M}$. In this case we call $\alpha$ a dynamical system on
${\mathcal M}$.
\end{definition}

We can define the infinitesimal generator of a dynamical system as
follows:

\begin{definition} Let $\alpha:\R\to U({\mathcal H})$ be a dynamical
system on ${\mathcal M}$, we define the infinitesimal generator
$\delta$ of $\alpha$ as a mapping $\delta:D(\delta)\subseteq
{\mathcal M}\to {\mathcal M}$ where
\[D(\delta)=\{{x\in {\mathcal M}:~ \lim_{t\to 0}\frac{\alpha_{t}(x)-x}{t}} ~
{\rm exists}\}\] and
\[\delta(x)=\lim_{t\to
0}\frac{\alpha_{t}(x)-x}{t},~ x\in D(\delta).\]
\end{definition}

Now we are ready to prove the main theorem of this paper:

\begin{theorem} Let ${\mathcal M}$ be a full Hilbert ${\mathcal A}$-module,
$\alpha $ be a dynamical system on ${\mathcal M}$ and $\delta$ be
the infinitesimal generator of $\alpha$. Then $D(\delta)$ is a
dense subspace of ${\mathcal M}$ and there exists a derivation
$d:D(d)\subseteq {\mathcal A}\to {\mathcal A}$ such that
$D(\delta)$ is a left $D(d)$-module and
$\delta(ax)=a\delta(x)+d(a)x, ~~a\in D(d), x\in D(\delta)$.
\end{theorem}

\begin{proof} By Hille-Yosida theorem \cite{H-P}, $D(\delta)$ is
a dense subspace of ${\mathcal M}$. Since $\alpha$ is a dynamical
system on ${\mathcal M}$ then for each $t\in \R$, the mapping
$\alpha_{t}:{\mathcal M}\to {\mathcal M}$ is a unitary. So there
exists $*$-isomorphism $\alpha_{t}^{'}:{\mathcal A}\to {\mathcal
A}$ such that $\langle \alpha_{t}(x),\alpha_{t}(y)\rangle
=\alpha_{t}^{'}(\langle x,y\rangle )$ and hence
$\alpha_{t}(ax)=\alpha_{t}^{'}(a)\alpha_{t}(x)~~(a\in {\mathcal
A}, x\in {\mathcal M}).$

Now we show that $\alpha{'}:\R\to Aut({\mathcal A})$ is a
$C^*$-dynamical system. For each $a\in {\mathcal A} , ~x\in
{\mathcal M} $ we have
$ax=\alpha_{0}(ax)=\alpha_{0}^{'}(a)\alpha_{0}(x)=\alpha_{0}^{'}(a)x$
thus by Lemma 2.3 $\alpha_{0}^{'}(a)=a$ for all $a \in {\mathcal
A}$. Therefore $\alpha_{0}^{'}=I$.

Also for all $t,s \in \R$ we have
\begin{eqnarray*}
\alpha_{t+s}^{'}(a)\alpha_{t+s}(x)&=&\alpha_{t+s}(ax)\\
&=&\alpha_{t}(\alpha_{s}(ax))\\
&=& \alpha_{t}(\alpha_{s}^{'}(a)\alpha_{s}(x))\\
&=&\alpha_{t}^{'}(\alpha_{s}^{'}(a))\alpha_{t}(\alpha_{s}(x))\\
&=&\alpha_{t}^{'}(\alpha_{s}^{'}(a))\alpha_{t+s}(x)
\end{eqnarray*}
and so $\alpha_{t+s}^{'}(a)=\alpha_{t}^{'}(\alpha_{s}^{'}(a))$.
Thus $\alpha_{t+s}^{'}=\alpha_{t}^{'}\alpha_{s}^{'}.$

Since for each $x\in {\mathcal M} $, $\alpha_{t}(x)\to x$ as $t\to
0$ we have
\[\|\alpha_{t}^{'}(a)x-ax\|\leq\|\alpha_{t}^{'}(a)x-\alpha_{t}^{'}(a)\alpha_{t}(x)\|
+\|\alpha_{t}^{'}(a)\alpha_{t}(x)-ax\|.\] Thus $\lim_{t\to 0}
\alpha_{t}^{'}(a)x=ax$ for all $x\in {\mathcal M}$ whence
$\lim_{t\to 0}\alpha_{t}^{'}(a)=a$ for all $a\in {\mathcal A}$.

Therefore $\alpha{'}:\R\to Aut({\mathcal A})$ is a $C^*$-dynamical
system on ${\mathcal A}$.

If $d$ is the infinitesimal generator of $\alpha{'}$ then for each
$a\in D(d), x\in D(\delta)$ we have
\begin{eqnarray*}
\lim_{t\to 0}\frac {\alpha_{t}(ax)-ax}{t}&=&\lim_{t\to 0}\frac
{a\alpha_{t}(x)-ax}{t}+
\lim_{t\to 0}\frac{\alpha{'}_{t}(a)\alpha_{t}(x)-a\alpha_{t}(x)}{t}\\
&=&a\lim_{t\to 0}\frac {\alpha_{t}(x)-x}{t}+ \lim_{t\to 0}\frac{\alpha{'}_{t}(a)-a}{t}\alpha_{t}(x)\\
&=&a\delta(x)+d(a)x
\end{eqnarray*}
So $ax\in D(\delta)$ , $\delta(ax)=a\delta(x)+d(a)x$. Furthermore
$D(\delta)$ is a left $D(d)$-module.
\end{proof}

\begin{example} Let ${\mathcal H}$ be a Hilbert space and $T\in B({\mathcal H})$ be
a self adjoint operator. By Stone's theorem there is a
one-parameter group of unitaries $t\mapsto \alpha_{t}
(\alpha_{t}\in B({\mathcal H}))$ which its infinitesimal generator
is $iT$. By Example 2.2, this is a dynamical system on Hilbert
$K({\mathcal H})$-module ${\mathcal H}$ with corresponding
dynamical system $t\mapsto \alpha_{t}{'}$ where
$\alpha_{t}{'}(V)=\alpha_{t}V\alpha_{t}^{*}$. Here
$\alpha_{t}^{*}$ denotes the adjoint of $\alpha_{t}$ in
$B({\mathcal H})$. If $\delta$ and $\delta^{'}$ are infinitesimal
generators of $\alpha$ and $\alpha{'}$ respectively then
$\delta=iT$ and by Theorem 4.3 we have
$\delta(Vx)=V\delta(x)+\delta^{'}(V)x$ and hence
$\delta(V(x))=V(\delta(x))+\delta^{'}(V)(x)$. Thus
\[\delta^{'}(V)(x)=\delta(V(x))-V(\delta(x))=iT(V(x))-iV(T(x)).\]
So \[\delta^{'}(V)=i(TV-VT)=i[T,V],~~ V\in D(\delta^{'}).\] The
density of $D(\delta^{'})$ in $K({\mathcal H})$ and the fact that
$\delta^{'}$ has no proper extension \cite{SAK} implies that
$\delta^{'}(V)=i[T,V],~~V\in K({\mathcal H})$.
\end{example}

We have shown that generalized derivations, like ordinary
derivations, behave well under some operations. In particular, 
closed derivations can be regarded as the infinitesimal generators
of some dynamical systems on Hilbert C*-modules. Since Hilbert
C*-modules naturally generalize C*-algebras, we may extend the
study of approximately dynamical systems to Hilbert C*-modules
\cite{NIK2}. We close the paper with a relevant question which may
be of special interest in the theory of approximately inner
dynamical systems \cite{SAK}:\\ Is a generalized derivation
approximately inner?
\\
{\bf Acknowledgement}: The authors would like to thank the referee for his/her useful comments and suggestions.

\end{document}